\documentclass{article}

\usepackage{arxiv}

\usepackage[utf8]{inputenc} 
\usepackage[T1]{fontenc}    
\usepackage{hyperref}       
\usepackage{url}            
\usepackage{booktabs}       
\usepackage{amsfonts}       
\usepackage{nicefrac}       
\usepackage{microtype}      
\usepackage{lipsum}		
\usepackage{amsthm}
\usepackage{graphicx}
\usepackage{doi}
\usepackage{physics}
\usepackage{amssymb}
\usepackage{pifont}

\newcommand{\pd}[2]{\frac{\partial #1}{\partial #2}}
\newcommand{\pdd}[2]{\frac{\partial^2 #1}{\partial #2^2}}
\newcommand{\bs}{\boldsymbol}
\newcommand{\mbf}{\mathbf}

\usepackage{setspace}
\title{The explicit constraint force method for optimal experimental design}

\date{} 					

\author{
    Conor Rowan \\
	Aerospace Engineering\\
	University of Colorado Boulder\\
    3775 Discovery Drive \\
	Boulder, CO 80309 \\
	\texttt{conor.rowan@colorado.edu} \\
}

\renewcommand{\headeright}{}
\renewcommand{\undertitle}{}
\renewcommand{\shorttitle}{}

\begin{document}
\maketitle

\begin{abstract}
    The explicit constraint force method (ECFM) was recently introduced as a novel formulation of the physics-informed solution reconstruction problem, and was subsequently extended to inverse problems. In both solution reconstruction and inverse problems, model parameters are estimated with the help of measurement data. In practice, experimentalists seek to design experiments such that the acquired data leads to the most robust recovery of the missing parameters in a subsequent inverse problem. While there are well-established techniques for designing experiments with standard approaches to the inverse problem, optimal experimental design (OED) has yet to be explored with the ECFM formulation. In this work, we investigate OED with a constraint force objective. First, we review traditional approaches to OED based on the Fisher information matrix, and propose an analogous formulation based on constraint forces. Next, we reflect on the different interpretations of the objective from standard and constraint force-based inverse problems. We then test our method on several example problems. These examples suggest that an experiment which is optimal in the sense of constraint forces tends to position measurements in the stiffest regions of the system. Because the responses---and thus the measurements---are small in these regions, this strategy is impractical in the presence of measurement noise and/or finite measurement precision. As such, our provisional conclusion is that ECFM is not a viable approach to OED.
\end{abstract}

\keywords{Explicit constraint force method \and Inverse problems \and Design of experiment \and Fisher information matrix}


\section{Introduction}
\label{intro}

\paragraph{} From Newton's law of gravitation to the Navier-Stokes and Maxwell equations, all scientific models depend on empirical parameters that are calibrated from measurement data. Unlike modern machine learning methods based on large neural networks, the models of engineering and physics make use of extensive prior knowledge, often in the form of conservation laws, which drastically reduces the number of empirical parameters to fit. In practice, ``inverse problems'' are solved for the model parameters, which philosopher and physicist Mario Bunge describes as identifying the ``hidden causes of observable effects'' \cite{bunge_inverse_2019}. It is interesting to note that the difference between problems viewed as data-fitting (through regression) and model calibration (with inverse problems) is largely a function of the number of parameters, though it is common to expect different generalization properties from these two types of models \cite{brunton_discovering_2016}. 

\paragraph{} Inverse problems can be probabilistic or deterministic in nature. For example, the Bayesian approach to inverse problems uses data to construct a probability distribution over the unknown parameters \cite{dashti_bayesian_2015}. By incorporating domain knowledge through a prior distribution on the parameters, the Bayesian approach is useful when measurement data is scarce. Furthermore, Bayesian methods are convenient for treating ill-posed problems, as multi-modal posterior distributions avoid overconfident point estimates of parameters which are non-unique \cite{sun_local_2022}. In contrast to the Bayesian approach, deterministic inverse problems seek parameters such that the predictions of the model agree as closely as possible with the measurement data \cite{vogel_computational_2002}. This approach has been particularly popular in the physics-informed machine learning literature in recent years, owing to its convenient incorporation into the standard pipeline of physics-informed neural network training. For example, the original physics-informed neural networks paper recovered the missing viscosity from the Navier-Stokes equations using flow data \cite{raissi_physics-informed_2019}.

\paragraph{} Solving an inverse problem relies on measurement data generated from an experiment. As the literature on identifiability shows, not all data provides sufficient information to recover the parameter(s) of interest \cite{miao_identifiability_2011, maclaren_what_2020, hong_global_2020}. Roughly speaking, the parameters of interest are identifiable if, with a given data set and parameterized model, there is a unique parameter setting that minimizes the error with the measurement data. For example, if two sets of parameters predict the same solution at the measurement locations, the problem is non-identifiable. The potential for non-identifiable inverse problems suggests that experimental efforts must be carried out with subsequent inverse problems in mind. There is not only freedom to choose the position of measurements, but also the size of the experimental coupon \cite{senig_investigation_2025}, as well as the position/intensity of source terms in the system \cite{krampe_optimized_2021}. Systematic approaches to design maximally informative experiments fall under the field of optimal experimental design (OED). As stated in \cite{huan_optimal_2024}, OED provides mathematical answers to the following questions:

\begin{quote}
    Where to place a sensor? What experimental conditions to impose? What quantity to observe? Do measurements need to be very precise, or would a noisier measurement suffice? When should measurements be made? And how much data should be collected? More broadly, what combination of observations would be most informative or useful—and how should we precisely define notions of ‘informative’ or ‘useful’ in the first place?
\end{quote}

In a number of fields, OED strategies have been successfully used to improve the quality of incoming data, and thus obtain more robust parameter estimates. In seismology, careful positioning of measurements with OED has been shown to reduce the computational cost of full-wave inversion problems \cite{mercier_designing_2025}. In solid mechanics, choosing optimal load paths for an elastoplastic sample leads to smaller uncertainties in the recovered constitutive model \cite{ricciardi_bayesian_2023}. Finally, OED techniques have been used to design field experiments in agricultural contexts, leading to more sustainable and profitable farms \cite{matavel_bayesian-optimized_2025}.

\paragraph{} While OED methods have been researched in a wide variety of application areas, they all rely on traditional formulations of the inverse problem. Recently, the explicit constraint force method (ECFM) was introduced as a novel approach to performing solution reconstruction in the presence of misspecified models \cite{rowan_physics-informed_2025}. ECFM introduces additional source terms into the system to enforce agreement between the parameterized model and the measurement data. The model parameters are then chosen such that the magnitude of the constraint force is minimal, and the value of the constraint force at the minimum provides a measure of consistency between the parameterized system and the data. This method was designed for solution reconstructions where the data is treated as ground truth, and the parameterized model is potentially mis-specified. Such situations arise with digital twins, where the underlying physics are often unknown or intractably complex \cite{tao_digital_2019}. Owing to the mathematical similarity between solution reconstruction and inverse problems, a follow-up work extended ECFM to inverse problems \cite{rowan_extending_2025}. Here, as in the case of solution reconstruction, the parameters are chosen such that the constraint force ensuring consistency with the measurement data is minimal. This contrasts with standard deterministic formulations, which seek parameters such that the mismatch between the system state and the data is minimal. 

\paragraph{} The purpose of this work is to explore OED with the constraint force formulation of this inverse problem. In Section 2, we review a traditional approach to experimental design based on the Fisher information matrix and its eigenvalues, and then, noting that ECFM is most natural in a deterministic setting, propose an interpretation of this problem which does not rely on probabilistic inference. Next, we review inverse problems for systems governed by partial differential equations (PDEs). In Section 3, we introduce the explicit constraint force method, and propose an extension of OED to the inverse problem based on constraint forces. In Section 4, we contrast the two perspectives on the inverse problem, arguing that ECFM provides parameters which are optimal in the sense of ``force,'' whereas the standard squared error approach is optimal in the sense of ``displacement.'' This section provides intuition which will be used in making sense of the results of the subsequent examples. In Section 5, we test the method on an idealized problem, for which we have analytical solutions. The results of the four variations of the inverse problem we study suggest that OED based on constraint forces favors measurement positions in the stiffest regions of the system. In Section 6, we discuss these results, and ultimately conclude that ECFM is not well-suited for OED.


\section{Optimal experimental design}
\label{oed}

\subsection{State of the art}

\paragraph{} The most common approaches to OED rely on the Fisher information matrix. To illustrate this, consider the following nonlinear regression problem:

\begin{equation}\label{regression}
    \mathbf{y} = \mathbf{f}(\mathbf{x}; \boldsymbol \epsilon) + \boldsymbol \eta,
\end{equation}

\noindent where $\mathbf{y}$ denotes the observed response of a system, $\mathbf f$ is a model parameterized by $\bs \epsilon$ which takes an input $\mathbf x$, and $\boldsymbol \eta$ is zero-mean measurement noise whose distribution is known. In general, this regression problem is nonlinear in the sense that the parameters $\boldsymbol \epsilon$ do not define a linear combination of basis functions. Suppose that we have $C$ observations $\{ \mathbf{x}_i , \mathbf{y}_i \}_{i=1}^C$. Using the distribution of the measurement noise, we can write the log likelihood as

\begin{equation}\label{loglike}
    \mathcal{L}(\boldsymbol \epsilon) = \sum_{i=1}^C \log p( \mathbf{y}_i |\mathbf{x}_i; \boldsymbol \epsilon),
\end{equation}

\noindent where $p(\mathbf{y}|\mathbf{x})$ is the density of the measurement noise. The Fisher information matrix $\mathbf I$ is defined as the expectation of the Hessian of the negative log likelihood \cite{pukelsheim_optimal_1993}:

\begin{equation*}
    \mathbf I(\boldsymbol{\epsilon}) = \mathbb{E}\qty( -\pdd{}{\boldsymbol \epsilon} \mathcal{L}(\boldsymbol \epsilon)).
\end{equation*}

Experimental design relies on the existence of ``experimental controls'' with which to adjust the setup of the experiment. These often correspond to measurement positions, but this is not a requirement. For example, an experimentalist might wish to determine the optimal size of a test sample given fixed measurement positions. Going forward, we denote the controls with a set of parameters $\boldsymbol \beta$. For the sake of exposition, we take the controls to be the measurement positions $\{ \mathbf{x}_i \}_{i=1}^C$. The goal of OED is then to determine measurement positions which are maximally informative given the parameterized model. While there are several distinct notions of ``informativeness,'' we choose to focus on the so-called E-criterion \cite{huan_optimal_2024}, which defines the OED problem as

\begin{equation}\label{e-crit}
    \underset{\boldsymbol \beta}{\text{argmax }} \text{min}(\text{eig}(\mathbf I(\boldsymbol \epsilon;\boldsymbol \beta)) ),
\end{equation}

\noindent where $\text{eig}( \bs \cdot)$ denotes the set of eigenvalues. The eigenvector corresponding to the minimum eigenvalue provides the direction of minimum curvature in the log likelihood loss landscape. As this curvature approaches zero, the parameters can vary in the corresponding eigenvector direction without changing the loss function, indicating non-identifiability. Thus, the E-criterion seeks to maximize the identifiability of the least identifiable direction on the loss surface. As it stands, there are two problems with Eq. \eqref{e-crit}: first, that the log likelihood depends on the observed response $\{ \mathbf y_i \}_{i=1}^C$, and second, that only the curvature of the loss at the optimal parameter setting $\boldsymbol \epsilon^*=\text{argmin } \mathcal L(\bs \epsilon)$ is of interest. Because the parameters are only optimal through the maximization of the likelihood, $\bs \epsilon^*$ is another quantity which requires knowledge of the data. If we take the experimental controls $\boldsymbol \beta$ to be the measurement positions, this means that in order to work with Eq. \eqref{e-crit}, we require $\mathbf{y}(\boldsymbol \beta)$. However, if the response of the system is known at any and all measurement points, there is no need to strategically position a discrete set of measurements with OED. These are fundamental tensions in OED---in order to \textit{truly} design an optimal experiment, the response of the system must be known in advance. But, if this condition is met, there is no need to design the experiment.

\paragraph{} In order to illustrate how these tensions are resolved, we work with a concrete example. We take the regression problem of Eq. \eqref{regression} and assume that the measurement noise is uncorrelated in the components of the system response and independent at each point in space:

\begin{equation*}
    \mathbf{y}_i = \mathbf{f}(\mathbf x_i ; \boldsymbol \epsilon) + \boldsymbol \eta_i, \quad \boldsymbol \eta_i \overset{i .i.d.}{\sim} N( \mathbf{0}, \sigma^2 \mathbf{I}),
\end{equation*}

\noindent where $\sigma^2$ is the variance in each coordinate direction. Given that the noise follows a multivariate normal distribution, we write

\begin{equation*}
    p(\mathbf{y}| \mathbf{x};\boldsymbol \epsilon) \propto \text{exp}\qty(-\frac{1}{2\sigma^2} \lVert \mathbf{y} - \mathbf f(\mathbf x ; \boldsymbol \epsilon)\rVert^2).
\end{equation*}

The log likelihood of Eq. \eqref{loglike} and its Hessian are straightforward to compute analytically. They are given by

\begin{equation}\label{normal}
    \mathcal{L}(\boldsymbol \epsilon) = -\frac{1}{2\sigma^2}\sum_{i=1}^C \lVert \mathbf y_i - \mathbf f( \mathbf x_i ; \boldsymbol \epsilon) \rVert^2, \quad \pdd{\mathcal{L}}{\boldsymbol \epsilon} = \frac{1}{\sigma^2} \sum_{i=1}^C \qty( -\pd{\mathbf f(\mathbf x_i)}{\boldsymbol \epsilon} \cdot \pd{\mathbf f(\mathbf x_i)}{\boldsymbol \epsilon} + (\mathbf y_i - \mathbf f( \mathbf x_i ; \boldsymbol \epsilon) ) \cdot \pdd{\mathbf f (\mathbf x_i)}{\boldsymbol \epsilon} ).
\end{equation}

Now, we note that to obtain the Fisher information matrix, we take an expectation of the negative Hessian in Eq. \eqref{normal}. Rewriting the model residual $(\mathbf y - \mathbf f)$ using the measurement noise, the Fisher information matrix is

\begin{equation}\label{fisher2}
    \mathbf I (\boldsymbol \epsilon) = \mathbb{E}\qty( \frac{1}{\sigma^2} \sum_{i=1}^C \qty( \pd{\mathbf f(\mathbf x_i)}{\boldsymbol \epsilon} \cdot \pd{\mathbf f(\mathbf x_i)}{\boldsymbol \epsilon} - \boldsymbol \eta_i \cdot\pdd{\mathbf f (\mathbf x_i)}{\boldsymbol \epsilon} )) = \frac{1}{\sigma^2} \sum_{i=1}^C \qty( \pd{\mathbf f(\mathbf x_i)}{\boldsymbol \epsilon} \cdot \pd{\mathbf f(\mathbf x_i)}{\boldsymbol \epsilon} ),
\end{equation}

\noindent which is a consequence of the deterministic structure of the regression model and the zero-mean property of the noise. Note that the data does not appear in Eq. \eqref{fisher2}, which lifts the requirement that it is known in advance. Note how the Fisher information matrix is constructed from $C$ rank 1 contributions, given the outer product structure of the two gradient terms. This means that its rank is at most $C$, and in the case that $|\boldsymbol \epsilon|>C$, the minimum eigenvalue will always be zero. However, the Fisher information matrix can also become rank deficient if $\partial \mathbf{f}(\mathbf{x}_i) / \partial \boldsymbol \epsilon \parallel \partial \mathbf{f}(\mathbf{x}_j) / \partial \boldsymbol \epsilon$ for $i\neq j$. As such, Eq. \eqref{fisher2} has the convenient property of discouraging coincident measurement points.

\paragraph{} Unless the model $\mathbf f$ is linear in $\boldsymbol \epsilon$, the model parameters appear in the Fisher information matrix. Recall that we care only about the curvature of the loss at the optimal parameter setting, but that the optimal parameters simultaneously require the data and are the goal of the entire OED process. Thus, we have not managed to avoid the issue of the unknown evaluation point in parameter space of the Fisher information matrix.  One straightforward fix to this problem is to endow the parameters with a prior distribution $\rho(\boldsymbol \epsilon)$. Instead of maximizing the minimum eigenvalue of the Fisher information at the solution, the E-criterion is modified to maximize the minimum eigenvalue of the expected Fisher information matrix:

\begin{equation}\label{expected_e-crit}
    \underset{\boldsymbol \beta}{\text{argmax }} \text{min}\qty(\text{eig}\qty( \int \mathbf I(\boldsymbol \epsilon;\boldsymbol \beta) \rho(\boldsymbol \epsilon) d\boldsymbol \epsilon) ).
\end{equation}

This strategy to perform OED for models which are nonlinear in the parameters is often called the Bayesian Fisher information matrix \cite{mentre_sparse-sampling_1995}. By way of example, in the case of normally distributed noise, a given prior on the parameters, and parameterized measurement positions, the E-criterion optimization problem with the Bayesian Fisher information matrix reads

\begin{equation*}
    \underset{ \boldsymbol \beta_1, \dots, \boldsymbol \beta_N}{\text{argmax }} \text{min}\qty(\text{eig}\qty( \int  \sum_{i=1}^N \qty( \pd{\mathbf f(\boldsymbol \beta_i)}{\boldsymbol \epsilon} \cdot \pd{\mathbf f(\boldsymbol \beta_i)}{\boldsymbol \epsilon} )\rho(\boldsymbol \epsilon) d\boldsymbol \epsilon) ).
\end{equation*}

\subsection{A deterministic perspective}

\paragraph{} We remark that Eq. \eqref{expected_e-crit} can be obtained without reference to probabilistic inference. Understanding this will be useful for developing an OED strategy relying on constraint forces, which are more natural in the deterministic setting. To see this, we solve the regression problem of Eq. \eqref{regression} with a standard mean-squared error objective:

\begin{equation}\label{Z}
    \mathcal{Z}(\boldsymbol \epsilon) = \frac{1}{2}\sum_{i=1}^C \lVert \mathbf y_i - \mathbf f( \mathbf x_i ; \boldsymbol \epsilon ) \rVert^2.
\end{equation}

As is well-known, the minimization of this objective is equivalent to maximum likelihood estimation for the parameters $\boldsymbol \epsilon$ with zero-mean Gaussian measurement noise. But, to forego the probabilistic perspective, we do not dwell on this connection. As before, let us assume that the experiment is parameterized with controls $\boldsymbol \beta$. A natural extension of the E-criterion to deterministic regression is given by the following optimization problem:

\begin{equation}\label{deterministic}
    \begin{aligned}
        \underset{\boldsymbol \beta, \boldsymbol \epsilon}{\text{argmax }} \text{min} \qty( \text{eig} \qty(  \pdd{\mathcal{Z}}{\boldsymbol \epsilon})) \\
        \text{s.t. } \pd{\mathcal Z}{\boldsymbol \epsilon} = \mathbf{0}.
    \end{aligned}
\end{equation}

This problem is interpreted as follows: find the experimental set-up $\boldsymbol \beta$ such that, when the error is minimized with respect to $\boldsymbol \epsilon$, the least identifiable direction is as identifiable as possible. It is interesting to note that, mathematically speaking, Eq. \eqref{deterministic} is equivalent to designing a structure against buckling. If $\mathcal{Z}$ is the total potential energy, $\boldsymbol \epsilon$ are the displacement degrees of freedom, and $\boldsymbol \beta$ are design parameters, this optimization problem attempts to avoid perturbations to an equilibrium configuration which do not generate stress \cite{liu_bayesian_2025}. Unfortunately, this connection is simply a curiosity, as the usual tactics for buckling analysis cannot be deployed on this problem. Unlike OED with the Fisher information matrix, there is no ambiguity as to what setting of $\boldsymbol \epsilon$ should be used in computing the Hessian matrix. This is handled by the constraint that the regression parameters are optimal for given experimental controls. However, we still suffer from not knowing the relationship between the measurements $\{ \mathbf y_i\}_{i=1}^C$ and the experimental controls $\boldsymbol \beta$. This is because both the constraint equation and Hessian involve the data:

\begin{equation}\label{dZ}
    \pd{\mathcal{Z}}{\boldsymbol \epsilon} = -\sum_{i=1}^C (\mathbf{y}_i - \mathbf f(\mathbf x_i ; \boldsymbol \epsilon)) \cdot \pd{\mathbf f(\mathbf x_i)}{\boldsymbol \epsilon}, \quad \pdd{\mathcal Z}{\boldsymbol \epsilon} = \sum_{i=1}^C\qty( \pd{\mathbf f(\mathbf x_i)}{\boldsymbol \epsilon} \cdot \pd{\mathbf f(\mathbf x_i)}{\boldsymbol \epsilon} - (\mathbf{y}_i - \mathbf f(\mathbf x_i ; \boldsymbol \epsilon)) \cdot \pdd{\mathbf f(\mathbf x_i)}{\boldsymbol \epsilon} ).
\end{equation}

There is no getting around the fact that we cannot know the data or optimal parameter setting in advance of doing the experiment. However, we can again assume some knowledge of the parameters through a prior $\rho(\boldsymbol \epsilon)$. If the prior is accurate, then the discrepancy between the model and the data will be small at all experimental settings $\boldsymbol \beta$. This means that the constraint $\partial \mathcal Z / \partial \boldsymbol \epsilon = \mathbf 0$ is approximately enforced. It also means that the second term in the Hessian of the objective of Eq. \eqref{dZ} is small. Neglecting this term yields the Gauss-Newton approximation of the Hessian, which is a common strategy in non-linear regression \cite{chen_hessian_2011}. Once again, we compute an expectation of a Hessian matrix over the prior on the model parameters. Taking the experimental parameters to be the measurement positions, this line of reasoning leads to the Bayesian Fisher information matrix from the probabilistic inference problem discussed above:

\begin{equation*}
    \underset{\boldsymbol \beta_1, \dots, \boldsymbol \beta _N}{\text{argmax }} \text{min}\qty(\text{eig}\qty( \int  \sum_{i=1}^C \qty( \pd{\mathbf f(\boldsymbol \beta_i)}{\boldsymbol \epsilon} \cdot \pd{\mathbf f(\boldsymbol \beta_i)}{\boldsymbol \epsilon} )\rho(\boldsymbol \epsilon) d\boldsymbol \epsilon) ).
\end{equation*}

Whereas the derivation of Fisher information matrix used an expectation and the assumption of zero mean noise to rid of the measurement data, we have asserted that a meaningful prior on the model parameters justifies the use of the Gauss-Newton approximation of the Hessian of the standard squared error regression objective. These arguments are rather similar, though moving out of the context of maximum likelihood estimation sets the stage for forthcoming discussions of OED with the explicit constraint force method. We remark that, regardless of how this matrix is obtained, it is positive semi-definite, which guarantees that its eigenvalues are non-negative. This guards against Eq. \eqref{fisher2} generating non-identifiable directions in the loss landscape by pushing a negative eigenvalue toward zero. We now discuss how this same deterministic line of reasoning applies to systems governed by PDEs.

\subsection{Designing experiments for PDE systems}

\paragraph{} The previous discussion on experimental design was restricted to nonlinear regression. Our goal now is to shift the discussion to inverse problems, which use measurement data to infer parameters of a physical system governed by ordinary or partial differential equations. Though inverse problems can certainly be thought of as nonlinear regression---where the model parameters define a regression manifold through the PDE solve---we believe the two problems are sufficiently different in practice to warrant their own discussions. In this work, we focus only on static PDEs. The system of interest is taken to obey the following abstract boundary value problem:

\begin{equation}\label{abstract_bvp}
\begin{aligned}
    \mathcal{N}( \mathbf{u})(\mathbf{x}) = \mathbf{0}, \quad \mathbf{x} \in \Omega,\\
    \mathcal{B}(\mathbf{u})(\mathbf{x}) = \mathbf{0}, \quad \mathbf{x} \in \partial \Omega,
\end{aligned}
\end{equation}

\noindent where $\mathbf{u} \in \mathbb{R}^{1,2,3}$ is the system state, $\mathbf{x}\in \mathbb{R}^{1,2,3}$ is the spatial coordinate, $\Omega$ is the domain on which the state is defined, $\mathcal{N}(\boldsymbol \cdot)$ is a differential operator, and $\mathcal{B}(\boldsymbol \cdot)$ is a boundary operator. The system state $\mathbf{u}(\mathbf{x})$ is taken to obey Eqs. \eqref{abstract_bvp} and is observed at $C$ distinct points. The measurements are given by 

\begin{equation}\label{measure}
    \mathbf{v}_i = \mathbf{u}(\mathbf{x_i}) + \boldsymbol{\eta}_i,
\end{equation}

\noindent where the measurement noise again has known distribution and is independent across observations. We collect all the measurements into a matrix with $\mathbf V = [\mathbf v_1 , \dots, \mathbf v_C]^T$. The analogue of the parametric function $\mathbf f$ is a parameterized boundary value problem given by 

\begin{equation}\label{parameterized_bvp}
\begin{aligned}
    \mathcal{G}( \mathbf{w} ; \boldsymbol \epsilon)(\mathbf{x}) = \mathbf{0}, \quad \mathbf{x} \in \Omega,\\
    \mathcal{Q}(\mathbf{w} ; \boldsymbol \epsilon)(\mathbf{x}) = \mathbf{0}, \quad \mathbf{x} \in \partial \Omega,
\end{aligned}
\end{equation}

\noindent where $\mathbf{w}$ is the state, $\mathcal{G}(\boldsymbol \cdot;\boldsymbol \epsilon)$ is the parameterized differential operator, and $\mathcal{Q}(\boldsymbol \cdot ; \boldsymbol \epsilon)$ is the parameterized boundary operator. To obtain a numerical solution, the state field is discretized with $\mathbf w(\mathbf x) \approx \mathbf{\hat w}(\mathbf{x}; \bs \theta)$ where $\bs \theta \in \mathbb R^N$ are the solution parameters. With the discretization of the solution, a numerical solution to Eq. \eqref{parameterized_bvp} is obtained either through minimizing the strong form error, solving the Galerkin weak form system, or, in the case that it exists, minimizing the variational energy. Regardless of how the numerical solution is approached, the solution is governed by a system of equations, which we denote as

\begin{equation}\label{res}
    \mathbf R( \bs \theta; \boldsymbol \epsilon) = \mathbf 0.
\end{equation}

When the discretization is linear, meaning that $\mathbf{\hat w} = \sum_{j=1}^{N} \theta_j \mathbf{f}_j(\mathbf{x})$ for some set of basis functions $\{ \mbf f_j(\mbf x)\}_{j=1}^{N}$, we can write a noiseless version of the measurement operation of Eq. \eqref{measure} as $\hat w_k(\mbf x_i) = \sum_{j=1}^N \theta_j f_{jk}(\mbf x_i)$, where $f_{jk}(\mbf x_i)$ is the $k$-th component of the $j$-th basis function evaluated at the $i$-th measurement position. This allows us to define the measurement operator $M_{kij}=f_{jk}(\mbf x_i)$. Analogous to the deterministic regression objective of Eq. \eqref{Z}, the inverse problem is

\begin{equation}\label{pde_inverse}
  \mathcal{Z}(\boldsymbol \epsilon) = \frac{1}{2} \lVert \mathbf M \bs \theta (\boldsymbol \epsilon) - \mathbf V \rVert^2 , \quad \underset{\boldsymbol \epsilon}{\text{argmin }} \mathcal{Z}(\boldsymbol \epsilon),
\end{equation}

\noindent where $\lVert \boldsymbol \cdot \rVert$ indicates the Frobenius norm in the case of vector-valued states, and thus matrix-valued data \cite{vogel_computational_2002}. Note that we think of the solution parameters $\bs \theta$ as a function of the model parameters $\boldsymbol \epsilon$ through the governing equation of Eq. \eqref{res}. When performing OED, we introduce experimental parameters $\boldsymbol \beta$ to build the experiment. These may appear in the measurement operator $\mathbf M$, or in the governing equation itself. See Table \ref{tab:parameters} to clarify the terminology we use for various types of parameters in the PDE-based experimental design problem. Using the E-criterion with the Gauss-Newton approximation of the Hessian and a prior on the model parameters, the OED problem is 

\begin{equation}\label{jinv}
    \underset{\boldsymbol \beta}{\text{argmax }} \text{min}\qty( \text{eig}\qty( \int  \mathbf{M} \pd{\bs \theta}{\boldsymbol \epsilon} \cdot \mathbf{M} \pd{\bs \theta}{\boldsymbol \epsilon}\rho(\boldsymbol \epsilon) d\boldsymbol \epsilon) ) =  \underset{\boldsymbol \beta}{\text{argmax }} \text{min}\qty( \text{eig}\qty( \mathbf{J}^{\text{INV}}(\boldsymbol \beta) ) ).
\end{equation}

The notation``$\text{INV}$'' is used to indicate that this OED problem arises from a standard formulation of the inverse problem. We note that sensitivity analysis through the governing equation of Eq. \eqref{res} can be used to obtain the derivatives $\partial \bs \theta / \partial \boldsymbol \epsilon$. This ensures that the solution parameters are always treated as an explicit function of the model parameters, i.e., $\bs \theta = \bs \theta(\boldsymbol \epsilon)$. The introduction of the measurement operator $\mathbf M $ and the necessity of computing $\bs \theta(\boldsymbol \epsilon)$ through a potentially non-linear solve are the primary differences between the nonlinear regression of the previous section and inverse problems for discretized PDE systems. In order to solve this optimization problem, we require the gradient of the minimum eigenvalue of the symmetric matrix $\mathbf{J}^{\text{INV}}$. Standard sensitivity analysis techniques for eigenproblems can be used in order to avoid the differentiating through a discontinuous $\text{min}(\boldsymbol \cdot)$ function. See Appendix \ref{sec: eigenvalues} for an exposition of eigenvalue sensitivity analysis.

\begin{table}[h]
    \centering
    \caption{Notation and meaning of the three types of parameters in the OED problem.}
    \label{tab:parameters}
    \begin{tabular}{lll}
        \toprule
        \textbf{Symbol} & \textbf{Name} & \textbf{Description} \\
        \midrule
        $\boldsymbol{\epsilon}$ & Model parameters & Material properties, source terms, boundary conditions, etc.\\
        $\boldsymbol \beta$ & Experimental parameters/controls & Measurement positions, sample size, point of load application, etc. \\
        $\bs \theta$ & Solution parameters & Coefficients in the discretized solution \\
        \bottomrule
    \end{tabular}
\end{table}

\section{Designing experiments with constraint forces}
\label{sec:cfs}

\paragraph{} The explicit constraint force method introduces source terms to the parameterized system of Eq. \eqref{parameterized_bvp} to enforce agreement between the state field and measurement data \cite{rowan_physics-informed_2025}. With ECFM, the parameterized boundary value problem is 

\begin{equation}\label{parameterized_bvp_ecfm}
\begin{aligned}
    \mathcal{G}( \mathbf{w} ; \boldsymbol \epsilon)(\mathbf{x}) + \sum_{i=1}^C \boldsymbol \lambda_i \Gamma (\mathbf{x} - \mathbf{x}_i) = \mathbf{0}, \quad \mathbf{x} \in \Omega,\\
    \mathcal{Q}(\mathbf{w} ; \boldsymbol \epsilon)(\mathbf{x}) = \mathbf{0}, \quad \mathbf{x} \in \partial \Omega,
\end{aligned}
\end{equation}

\noindent where $\Gamma(\mathbf{x})$ is a user-specified form of the constraint force, $\{ \mathbf{x}_i \}_{i=1}^C$ are the measurement locations, and $ \{\boldsymbol \lambda_i\}_{i=1}^C$ scale the constraint forces. The discretized form of Eq. \eqref{parameterized_bvp_ecfm} is given by 

\begin{equation*}
    \mathbf{R}( \bs \theta ; \boldsymbol \epsilon) + \boldsymbol \Gamma \boldsymbol \lambda = \mathbf{0},
\end{equation*}

\noindent where $\boldsymbol \Gamma$ stores the force vectors corresponding to the constraint forces at each of the measurement locations. As shown in \cite{rowan_extending_2025}, the ECFM approach to inverse problems is 

\begin{equation}\label{ecfm}
    \begin{aligned}
        \underset{\bs \theta, \boldsymbol \epsilon,\boldsymbol \lambda}{\text{argmin }} \frac{1}{2} \lVert \boldsymbol \lambda \rVert^2 \\
        \text{s.t. } \mathbf{ R}( \bs \theta; \boldsymbol \epsilon) +  \boldsymbol \Gamma \boldsymbol \lambda = \mathbf{0}, \quad \mathbf M  \bs \theta - \mathbf V = \mathbf 0.
    \end{aligned}
\end{equation}

In words, this optimization problem finds the model parameters $\boldsymbol \epsilon$ such that the least amount of additional ``forcing'' is required to make the solution parameters $\bs \theta$ consistent with the measurement data $\mathbf V$. We say that the parameterized model of Eq. \eqref{parameterized_bvp_ecfm} is ``consistent'' with the true system when the solution to Eq. \eqref{ecfm} has a constraint force of $\boldsymbol \lambda = \mathbf 0$. This means that there is at least one setting of the model parameters that exactly reproduces the data. Admittedly, the ECFM formulation of the inverse problem appears more complex than the standard approach of Eq. \eqref{pde_inverse}. When the parameterized model is consistent with the measurement data, these two methods will yield the same estimate of the solution and model parameters. However, differences between the two approaches arise when the parameterized model cannot reproduce the measurement data.

\paragraph{} As before, we assume for the sake of exposition that the experimental controls determine measurement positions. This means that $\bs \beta$ shows up in the measurement operator and the constraint force vectors. The analogous OED problem to Eq. \eqref{deterministic} is

\begin{equation}\label{ecfm_oed}
    \begin{aligned}
        \underset{\bs \theta, \boldsymbol \epsilon, \boldsymbol \lambda, \boldsymbol \beta}{\text{argmax }} \text{min} \qty( \text{eig}\qty ( \pdd{}{\boldsymbol \epsilon} \frac{1}{2} \lVert \boldsymbol \lambda \rVert^2) ) \\
        \text{s.t. } \pd{}{\boldsymbol \epsilon} \frac{1}{2} \lVert \boldsymbol \lambda \rVert^2 = \mathbf 0, \quad  \mathbf{ R}( \boldsymbol \theta; \boldsymbol \epsilon) +  \boldsymbol \Gamma(\bs \beta) \boldsymbol \lambda = \mathbf{0}, \quad \mathbf M(\bs \beta)  \bs \theta - \mathbf V(\bs \beta) = \mathbf 0.
    \end{aligned}
\end{equation}

This optimization problem is interpreted as follows: find the experimental controls such that perturbations to the model parameters around the minimum constraint force solution in the direction of minimal change cause the largest possible change to the constraint force magnitude. Like the standard inverse problem, this objective works to maximize the identifiability of the least identifiable direction in $\boldsymbol \epsilon$ space, though the notion of identifiability is distinct. This is a very complex optimization problem. And, just like OED with the standard inverse problem, we require the dependence of the data on the experimental controls in order to truly solve it. Thus, we need to somehow modify Eq. \eqref{ecfm_oed} in order to eliminate the data dependence. As before, we can do this with the help of a prior over the model parameters. To simplify exposition, we take the physics of the parameterized governing equation to be linear. The ECFM governing equation is 

\begin{equation*}
    \mbf K(\bs \epsilon) \bs \theta + \bs \Gamma(\bs \beta) \bs \lambda - \mbf F(\bs \epsilon) = \mbf 0,
\end{equation*}

\noindent where $\mbf F(\bs \epsilon)$ represents the force vector corresponding to a source term in the system and $\mbf K(\bs \epsilon)$ is the discretized differential operator. Both the source and differential operator are taken to depend on the model parameters $\bs \epsilon$ in some unspecified way. Now, we define an approximation of the data generating process using the prior over model parameters as

\begin{equation*}
    \mbf {\tilde V}(\bs \beta) = \mbf M(\bs \beta) \int \mbf K(\bs \epsilon)^{-1} \mbf F(\bs \epsilon) \rho(\bs \epsilon) d\bs \epsilon.
\end{equation*}

This is the average response of the parameterized system at the measurement positions defined by $\bs \beta$. For this approximation to be accurate, the parameterized model must be able to accurately reproduce the measurement data, and the prior over model parameters must be approximately centered around their true value. With a model of the data generation process, we might now find the optimal model parameters $\bs \epsilon^*(\bs \beta)$ as a function of the experimental controls by minimizing the constraint force magnitude. However, given the approximation made in the data generation, it seems needless to enforce this level of granularity in the point around which the Hessian in Eq. \eqref{ecfm_oed} is computed. Thus, we again use the prior to rid the OED optimization problem of a constraint, namely that the constraint is minimized. This is analogous to removing the constraint that the squared error is minimal in the previous section. An approximate solution to Eq. \eqref{ecfm_oed} is then provided by 

\begin{equation}\label{ecfm_oed_approx}
    \begin{aligned}
        \underset{\bs \theta, \boldsymbol \lambda, \boldsymbol \beta}{\text{argmax }} \text{min} \qty( \text{eig}\qty ( \int \pdd{}{\boldsymbol \epsilon} \frac{1}{2} \lVert \boldsymbol \lambda \rVert^2 \rho(\bs \epsilon) d \bs \epsilon) ) \\
        \text{s.t. }  \quad  \mathbf{ K}(\bs \epsilon) \boldsymbol \theta +  \boldsymbol \Gamma(\bs \beta) \boldsymbol \lambda - \mbf F(\bs \epsilon)= \mathbf{0}, \quad \mathbf M(\bs \beta)  \bs \theta - \mathbf {\tilde V}(\bs \beta) = \mathbf 0.
    \end{aligned}
\end{equation}

Unlike the OED approach based on the standard inverse problem, there is no clear way to make the dependence of the constraint force objective on the measurement data disappear. We can remove the constraint of this optimization problem as follows. The constraint force and solution parameters are governed by 

\begin{equation}\label{cfgovern}
    \begin{bmatrix}
        \mbf K( \bs \epsilon) &  \bs \Gamma(\bs \beta) \\ \mbf M(\bs \beta) & \mbf 0 
    \end{bmatrix} \begin{bmatrix}
        \bs \theta \\ \bs \lambda
    \end{bmatrix} - \begin{bmatrix}
        \mbf F(\bs \epsilon) \\ \mbf {\tilde V}(\bs \beta)
    \end{bmatrix} = \mbf 0.
\end{equation}

This means that the constraint force can be written as an explicit function of the model parameters and experimental controls as

\begin{equation*}
    \bs \lambda = \mbf T \begin{bmatrix}
        \mbf K( \bs \epsilon) &  \bs \Gamma(\bs \beta) \\ \mbf M(\bs \beta) & \mbf 0 
    \end{bmatrix}^{-1}\begin{bmatrix}
        \mbf F(\bs \epsilon) \\ \mbf {\tilde V}(\bs \beta)
    \end{bmatrix},
\end{equation*}

\noindent where $\mbf T$ is a truncation operator which picks out the constraint forces from the solution vector $[ \bs \theta, \bs \lambda]^T$. The final form of the OED optimization problem with constraint forces is

\begin{equation}\label{ecfm_oed_approx_final}
    \underset{  \boldsymbol \beta}{\text{argmax }} \text{min} \qty( \text{eig}\qty ( \int \pdd{}{\boldsymbol \epsilon} \frac{1}{2} \lVert \boldsymbol \lambda(\bs \epsilon, \bs \beta) \rVert^2 \rho(\bs \epsilon) d \bs \epsilon) ) .
\end{equation}

Eq. \eqref{ecfm_oed_approx_final} should be very exciting to enthusiasts of sensitivity analysis. We require the gradient of the minimum eigenvalue of the expectation of a Hessian of a quantity computed through a constrained PDE solve. As shown in Appendix \ref{sec: ecfm_sens}, this requires third sensitivity derivatives of the governing equation for the constraint force. In the following section, we take a step back and interpret physically the two formulations of the inverse problem we have discussed. This intuition will assist in making sense of the results of the example problems.


\section{Two conceptions of the inverse problem}
\label{philosophy}

\paragraph{} We claim that Eq. \eqref{pde_inverse} furnishes model parameters that are optimal in the sense of displacement, whereas the ECFM objective of Eq. \eqref{ecfm} finds model parameters that are optimal in the sense of force. What is meant by this? Abstractly, we believe that the solution to an inverse problem should be conceptualized as model parameters which minimize the discrepancy of a parameterized system with measurement data. There is, however, no single way to measure discrepancy. In the case of Eq. \eqref{pde_inverse}, the model parameters $\boldsymbol \epsilon$ are chosen such that the system state minimizes the error with the measurement data. To use an analogy with elasticity, we call the system state the ``displacement,'' but recognize that in general it could be another quantity such as temperature, velocity, etc. As a thought experiment, imagine an experimentalist adjusting the model parameters with knobs until the equilibrium configuration of an elastic system minimizes the distance to the measured displacements. When a setting is found such that small adjustments to each of the knobs increase the sum of squared errors with the measurement data, a solution is obtained. Here, a parameterized system is consistent with the measurement data when the system exactly reproduces the data at the optimal setting of the model parameters.

\begin{figure}[hbt!]
\centering
\includegraphics[width=.95\textwidth]{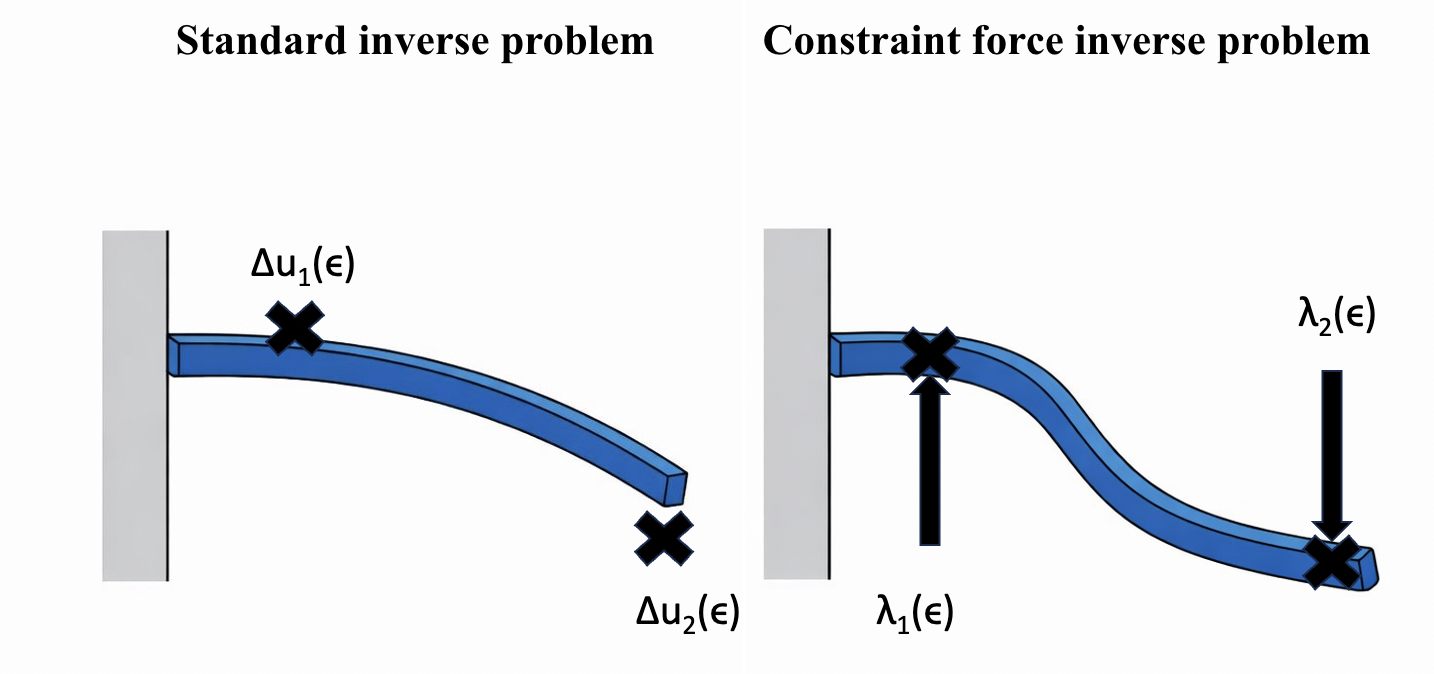}
\caption{In the case of the standard inverse problem, the model parameters $\boldsymbol \epsilon$ are chosen to minimize $\sum_i\Delta u_i^2(\bs \epsilon)$, which is the displacement discrepancy. With ECFM, the model parameters minimize the magnitude of the constraint force, given by $\sum_i \lambda_i^2(\bs \epsilon)$.}
\label{stacking}
\end{figure}

\paragraph{} ECFM provides another measure of discrepancy between the parameterized system and the data. Instead of minimizing the error in the displacements, ECFM minimizes the magnitude of the force required to make the system pass through the measurement data. In this case, the experimentalist adjusts the model parameter knobs, and at each setting, stacks weights on the structure until it agrees with the measurement data. When small adjustments to the knobs do not change the total amount of weight required to enforce agreement, a solution is obtained. If the parameterized system is consistent, there is a parameter setting for which the constraint force is zero, as there is no discrepancy in the system state to correct. Minimizing these two measures of discrepancy provides the same model parameters only when the system is consistent. In the case of inconsistency, the constraint force objective prioritizes discrepancies in stiffer regions of the system, as these require larger constraint forces to correct. The standard inverse problem treats all displacement errors equally, regardless of their location in the domain. Upon reflection, this egalitarian attitude is somewhat counterintuitive: a displacement discrepancy near the root of a clamped beam is weighted just as heavily as a discrepancy near the tip. In the case of noiseless measurements, minor failures to reproduce the data at the root may indicate model mis-specification to a greater extent than comparatively large errors near the tip. Standard formulations of the inverse problem---whose notion of consistency revolves around displacement rather than force---do not provide feedback of this sort. On the other hand, ECFM provides a more straightforward measure of the mis-specification, as the magnitude of the constraint force objective is sensitive to both the magnitude of the displacement discrepancies and their location. See Figure \ref{stacking} for a visualization of the two different formulations of the inverse problem. In the opinion of the author, the differences between the standard inverse problem and ECFM are largely philosophical, as there is no principled reason to believe that inverse problems in the case of mis-specified models provide meaningful parameter estimates. The two methods provide different routes to the same solution when the model is consistent with the data.

\section{Testing the method}
\label{sec:examples}

\paragraph{} We work with a simple model problem to test the proposed method for experimental design with constraint forces. This allows us to work entirely with analytic solutions, which provide clear insight into the behavior of the constraint force method. The model problem is given by

\begin{equation*}
    k\pdd{u}{x} + b = 0, \quad u(0)=0, \quad k\pd{u}{x}(1)=p.
\end{equation*}

The solution to this boundary value problem is

\begin{equation*}
    u(x) = -\frac{b}{2k}x^2 + \frac{p+b}{k}x.
\end{equation*}

In these examples, we assume that we have full access to the data generation process. For the reasons discussed above, this is an unrealistic assumption. Indeed, there would be no reason to optimize the design of the experiment if measurements of the system could be known in advance. We work with this idealized problem in order to expose the behavior of the constraint force method, not to solve a realistic problem. Additionally, we assume that the measurement is taken at a single point $\beta \in [0,1]$. As such, the data is given by 

\begin{equation*}
    v(\beta) = u(\beta) = -\frac{b}{2k}\beta^2 + \frac{p+b}{k}\beta.
\end{equation*}

In the following problems, we assume that there is a single model parameter to estimate. We remark that in the case of a single measurement, single model parameter, and analytic solution for the displacement, the E-criterion with the Fisher information matrix is simply

\begin{equation}\label{simple_jinv}
    \underset{\beta}{\text{argmax }} \int \qty( \pd{w(\beta)}{\epsilon} )^2 \rho( \epsilon) d  \epsilon.
\end{equation}

\noindent where $w(\beta)$ is the prediction of the parameterized model and $\epsilon$ is the scalar model parameter to be estimated. This amounts to positioning the measurement where the prediction is most sensitive to the model parameter on average. Now, we experiment with the constraint force formulation of the OED problem with different components of the model problem parameterized.

\subsection{Parameterized boundary condition}

\paragraph{} In the first test case, we take the right boundary condition to be the quantity of interest in the inverse problem. Per the ECFM approach to inverse problems, we introduce a constraint force centered at the measurement position. The governing equation is

\begin{equation*}
    k \pdd{w}{x} + b + \lambda \delta(x-\beta)= 0, \quad w(0)=0, \quad k\pd{w}{x}(1)=\epsilon.
\end{equation*}

Going forward, we choose to work with impulse constraint forces. This is a straightforward way to allow the measurement to be positioned anywhere in the domain without adjusting the width of the constraint force. The analytic solution to the governing equation is 

\begin{equation*}
    w(x) = -\frac{b}{2k}x^2 + \frac{\lambda}{k}( x - \text{max}(0,x-\beta)) + \frac{\epsilon+b}{k}x.
\end{equation*}

The constraint force $\lambda$ is determined as a function of $\epsilon$ and $\beta$ by setting $w(\beta)=v(\beta)$. In this case, we find that $\lambda = p - \epsilon $. With the single measurement and single model parameter, Eq. \eqref{ecfm_oed_approx_final} becomes

\begin{equation*}
    \underset{\beta}{\text{argmax }} \int \frac{1}{2}\pdd{}{\epsilon} (p-\epsilon)^2 \rho(\epsilon) d\epsilon.
\end{equation*}

The constraint force method is indifferent to where the measurement is positioned, given that $\beta$ does not appear in the expression for the constraint force magnitude. For the standard inverse problem, we have that

\begin{equation*}
    \pd{w(\beta)}{\epsilon} = \frac{\beta}{k},
\end{equation*}

\noindent which, in combination with Eq. \eqref{simple_jinv}, defines the experiment to be optimal when $\beta=1$. This is simply the position at which the model prediction changes the most with the applied traction. In contrast, the constraint force objective states that the same effort is required to correct the discrepancy with the data at all points in the system. Thus, the curvature of the constraint force objective is equivalent at all measurement positions. Practically speaking, it is counterintuitive that a measurement near the zero Dirichlet boundary condition on the left is as informative as a measurement on the right end where the displacement is much larger. Implicitly, this result is relying on arbitrary precision of the measuring device and the assumption of zero noise. Obviously, these conditions are not met in practice.

\subsection{Parameterized source term}

\paragraph{} In this second example, we take the boundary conditions to be known, but the distributed load to be estimated in the inverse problem. The governing equation is 

\begin{equation*}
    k \pdd{w}{x} + \epsilon + \lambda \delta(x-\beta)= 0, \quad w(0)=0, \quad k\pd{w}{x}(1)=p.
\end{equation*}

The exact solution is given by 

\begin{equation*}
    w(x) = -\frac{\epsilon}{2k}x^2 + \frac{\lambda}{k}( x - \text{max}(0,x-\beta)) + \frac{p+\epsilon}{k}x.
\end{equation*}

The constraint force is determined as a function of the model parameter and the measurement position by setting $w(\beta) = v(\beta)$. The constraint force is 

\begin{equation*}
    \lambda = \frac{\beta}{2}(\epsilon-b)+b-\epsilon.
\end{equation*}

The OED problem based on constraint forces is 

\begin{equation*}
    \underset{\beta}{\text{argmax }} \int \qty[\qty(\pd{\lambda}{\epsilon})^2 + \lambda \pdd{\lambda}{\epsilon}] \rho(\epsilon)d \epsilon = \underset{\beta}{\text{argmax }} \int \qty( \frac{\beta}{2}-1)^2 \rho(\epsilon)d \epsilon = \underset{\beta}{\text{argmax }} \qty( \frac{\beta}{2}-1)^2.
\end{equation*}

Regardless of the prior distribution, the solution to the optimization problem is obtained with $\beta=0$. The constraint force method determines the experiment to be optimal when the measurement coincides with the Dirichlet boundary. The Fisher information matrix will again position the measurement at the right end of the domain, as this is where the solution is most sensitive to the model parameter. The two OED methods make opposite determinations of optimality. A measurement at the left boundary does not make sense in practice, suggesting that the constraint force objective is problematic for OED.

\subsection{Parameterized material}

\paragraph{} Now, we take the material property to be the parameter of interest in the inverse problem. The governing equation is 

\begin{equation*}
    \epsilon \pdd{w}{x} + b + \lambda \delta(x-\beta)= 0, \quad w(0)=0, \quad \epsilon\pd{w}{x}(1)=p.
\end{equation*}

The exact solution is given by

\begin{equation*}
    w(x) = - \frac{b}{2\epsilon}x^2 +  \frac{\lambda}{\epsilon}( x - \text{max}(0,x-\beta)) + \frac{p+b}{\epsilon}x.
\end{equation*}

The constraint force is determined as a function of the model parameter and measurement position by setting $w(\beta)=v(\beta)$:

\begin{equation*}
    \lambda = \qty( \frac{b\beta}{2} - p - b )\qty( 1 - \frac{\epsilon}{k}).
\end{equation*}

Given the linear dependence of the constraint force on the model parameter, the constraint force OED objective is

\begin{equation*}
    \underset{\beta}{\text{argmax }} \int \qty(\pd{\lambda}{\epsilon})^2\rho(\epsilon) d\epsilon = \underset{\beta}{\text{argmax }} \frac{1}{k^2}\qty( p + b \qty(1 - \frac{\beta}{2}) )^2.
\end{equation*}

It is no longer possible to find the measurement position which maximizes this quadratic expression by inspection. Note that we require $0<\beta<1$, meaning this is a constrained optimization problem. Dropping the constant $1/k^2$, which does not affect the maximizer, we expand the polynomial objective as

\begin{equation*}
    \qty( p + b \qty(1 - \frac{\beta}{2}) )^2 = \frac{b^2}{4}\beta^2 - (pb+b^2)\beta + (p+b)^2.
\end{equation*}

Note that this function is convex in $\beta$, as $b^2/4>0$. By taking the derivative, the minimum is located at 

\begin{equation*}
    \beta^* = 2\qty(\frac{p}{b}+1).
\end{equation*}

Given that $\beta \in [0,1]$, and using convexity, we can see that if $\beta^*>1/2$, the maximum of the constraint force curvature is obtained at $\beta=0$. If $\beta<1/2$, the maximum is obtained at $\beta=1$. We cannot determine the measurement position without knowledge of the values of the traction and source term. Note that the optimal measurement position does not depend continuously on these values, it jumps suddenly between the left and right ends of the domain.

\subsection{Mis-specified model}

\paragraph{} In this example, the parameterized model is inconsistent with the true data-generating model. The source term is parameterized, and we falsely assume a homogeneous Neumann boundary condition on the right. The governing equation is then

\begin{equation*}
    k \pdd{w}{x} + \epsilon + \lambda \delta(x-\beta)= 0, \quad w(0)=0, \quad k\pd{w}{x}(1)=0.
\end{equation*}

The exact solution is given by 

\begin{equation*}
    w(x) = -\frac{\epsilon}{2k}x^2 + \frac{\lambda}{k}( x - \text{max}(0,x-\beta)) + \frac{\epsilon}{k}x.
\end{equation*}

By setting $w(\beta)=v(\beta)$, the constraint force is determined to be 

\begin{equation*}
    \lambda = \frac{\beta}{2}(\epsilon - \beta ) + p + b - \epsilon.
\end{equation*}

The OED problem based on constraint forces is 

\begin{equation*}
    \underset{\beta}{\text{argmax }} \int \qty[\qty(\pd{\lambda}{\epsilon})^2 + \lambda \pdd{\lambda}{\epsilon}] \rho(\epsilon)d \epsilon = \underset{\beta}{\text{argmax }} \int \qty( \frac{\beta}{2}-1)^2 \rho(\epsilon)d \epsilon = \underset{\beta}{\text{argmax }} \qty( \frac{\beta}{2}-1)^2.
\end{equation*}

Just as in the case of the parameterized source term when the model was correctly specified, the optimal measurement position is $\beta=0$. The constraint force-based OED problem once again pushes the measurement to the stiffest part of the domain.


\section{Discussion}
\label{sec:discussion}

\paragraph{} After reviewing the standard approach to OED with the Fisher information matrix, we showed how the E-criterion can be derived from a deterministic inverse problem using the Gauss-Newton approximation of the Hessian. In the deterministic setting, the E-criterion is understood to design experiments which lead to well-conditioned inverse problems. The deterministic perspective allowed us to extend the E-criterion to the ECFM formulation of the inverse problem. Given that an optimization problem for the minimum constraint force governs the model parameters with ECFM, we proposed an OED objective that maximized the minimum eigenvalue of the Hessian of the constraint force magnitude. We discussed how a prior over model parameters could be used to eliminate the OED problem's dependence on the unknown experimental data. Before proceeding to test the proposed method, we clarified how the two approaches to the inverse problem differ in their measures of the system's discrepancy with the data. We then tested the ECFM approach to OED on an idealized problem, in which we assumed that the data generating process was specified in advance of performing the experiment. Though unrealistic, this allowed us to study the fundamental behavior of the method, free of approximations. We found that optimal experiments in the ECFM sense rely on measurements at the stiffest position of the system. In the context of the example problem, this meant positioning measurements coincident with the Dirichlet boundary. This contrasts with the standard OED approach based on the Fisher information matrix, which positioned measurements where the system state was most sensitive to the model parameter. Obviously, taking measurements at a Dirichlet boundary is nonsensical---our proposed method does this because the constraint force magnitude is most sensitive to the model parameter at this ultra-stiff position in the system. The flaw in our approach is that the constraint force is not a quantity which is actually measured, rather it is a construct to measure consistency of the parameterized system with the data. However, even if the constraint force were measured---meaning that a real source term was introduced to a system in order for it to respect the measurement data---measurements in stiff regions would require arbitrarily high precision to determine the necessary constraint force. With these considerations in mind, we believe the idealized problem we have studied to show that ECFM is not a viable approach to designing optimal experiments. 

\paragraph{} While \cite{rowan_extending_2025} shows that ECFM can be used successfully to solve inverse problems, we fear that our conclusions here are a death knell for the ECFM formulation of inverse problems. How good could an approach be if, when ``optimized,'' it produces nonsense? More practically, where should measurements \textit{really} be positioned for an ECFM inverse problem? Perhaps the E-criterion, though intuitive as a measure of the ease of the optimization problem, lead us astray from the very beginning. Future work will likely halt for all but the most zealous ECFM devotees. Such stubborn and tragic individuals will focus their efforts on devising more exotic means to perfect the constraint force minimization problem. And, while the rest of us hum along, hemmed in by walls we fear but no longer see, the future of the constraint force is being authored by these chosen few. Whether it be a revolution, or but a hope, quickly dispersed by time's inexorable flow---this remains to be seen. In the meantime, to those among us who still have the gall to dream, I offer a prayer: \textit{the tyrannical reign of displacement will come to an end. Long live the constraint force!}

\appendix

\counterwithin*{equation}{section}
\renewcommand\theequation{\thesection\arabic{equation}}

\section{Eigenvalue sensitivity analysis}
\label{sec: eigenvalues}

\paragraph{} This appendix demonstrates how to compute the gradient of the minimum eigenvalue of a symmetric matrix $\mbf J$ with respect to experimental controls $\bs \beta$. We can accomplish this by writing the eigen-equation corresponding to the minimum eigenvalue, and differentiating. Calling the minimum eigenvalue $\mu$ and the corresponding eigenvector $\mathbf q$, this reads

\begin{equation*}
    \mathbf 0 = \pd{}{\boldsymbol \beta}( \mathbf{J} \mathbf q - \mu \mathbf q) = \pd{\mathbf J} {\boldsymbol \beta} \mathbf q + \mathbf J \pd{\mathbf q}{\boldsymbol \beta} - \pd{\mu}{\boldsymbol \beta} \mathbf q - \mu \pd{\mathbf q}{\boldsymbol \beta}.
\end{equation*}

Multiplying by $\mathbf q^T$ on the left and using orthonormality of the eigenvectors, the sensitivity derivative for the minimum eigenvalue is given by 

\begin{equation*}
    \pd{\mu}{\boldsymbol \beta} = \mathbf q^T \pd{\mathbf J}{\boldsymbol 
    \beta} \mathbf q.
\end{equation*}

This result holds provided the minimal eigenvalue does not collide or change ordering with the neighboring eigenvalue. Under this condition, the sensitivity of $\mu$ is well-defined.

\section{ECFM sensitivity analysis}
\label{sec: ecfm_sens}

\paragraph{} The following example is taken from \textit{Zen and the Art of Sensitivity Analysis} \cite{rowan_okay_2025}. We are interested in  computing the gradient of the minimum eigenvalue of the expected Hessian of the constraint force magnitude. Recall that the constraint force is treated as an explicit function of the model parameters $\bs \epsilon$ and the experimental controls $\bs \beta$. The matrix of interest is 

\begin{equation*}
    \mbf J =  \int \pdd{}{\boldsymbol \epsilon} \frac{1}{2} \lVert \boldsymbol \lambda(\bs \epsilon, \bs \beta) \rVert^2 \rho(\bs \epsilon) d \bs \epsilon.
\end{equation*}

Appendix \ref{sec: eigenvalues} shows the crux of sensitivity analysis for the minimum eigenvalue is obtaining $\partial \mbf J / \partial \bs \beta$. Things are about to get hairy---to keep things slightly more readable, we introduce the following convention. Indices corresponding to the model parameters are Greek, indices corresponding to the experimental controls are capitalized Latin, and indices corresponding to the solution parameters and constraint forces are lowercase Latin. The gradient of the expected constraint force magnitude is

\begin{equation}\label{JB}
\begin{aligned}
    \pd{J_{\alpha \gamma}}{\beta_K} = \int \pd{}{\beta_K}\qty( \pd{\lambda_{\ell}}{\epsilon_{\alpha} }\pd{\lambda_{\ell}}{\epsilon_{\gamma}} + \lambda_{\ell} \frac{\partial^2 \lambda_{\ell}}{\partial \epsilon_{\alpha} \partial \epsilon_{\gamma}}) \rho(\bs \epsilon) d\bs \epsilon = \\\int \qty(\frac{\partial^2 \lambda_{\ell}}{\partial \epsilon_{\alpha} \partial \beta_K} \pd{\lambda_{\ell}}{\epsilon_{\gamma}} + \pd{\lambda_{\ell}}{\epsilon_{\alpha}} \frac{\partial^2 \lambda_{\ell}}{\partial \epsilon_{\gamma} \partial \beta_K} + \pd{\lambda_{\ell}}{\beta_K} \frac{\partial^2 \lambda_{\ell}}{\partial \epsilon_{\alpha} \partial \epsilon_{\gamma}} + \lambda_{\ell} \frac{\partial^3 \lambda_{\ell}}{ \partial \epsilon_{\alpha} \partial \epsilon_{\gamma} \partial \beta_K}) \rho(\bs \epsilon) d \bs \epsilon .
\end{aligned}
\end{equation}

To compute the gradient in Eq. \eqref{JB}, we require (i) the constraint force, (ii) the first sensitivity derivative of the constraint force with respect to the model parameters, (iii) the first sensitivity derivative of the constraint force with respect to the experimental controls, (iv) the second sensitivity derivative of the constraint force with respect to the model parameters, (v) the mixed sensitivity derivative of the constraint force with respect to the model parameters and the controls, and (vi) the sensitivity derivative with respect to the experimental controls of the second sensitivity derivative of the constraint force with respect to the model parameters. Suppose we write the governing equation of Eq. \eqref{cfgovern} as $D_{ij}(\bs \epsilon, \bs \beta)y_j - Q_i(\bs \epsilon, \bs \beta) = 0$, where the solution vector is $\mbf y=[\bs \theta, \bs \lambda]^T$. The constraint and solution parameters are governed by

\begin{equation*}
    D_{ij}y_j - Q_i = 0.
\end{equation*}

Using this solution, the sensitivity derivative (ii) can be computed with

\begin{equation*}
    \pd{D_{ij}}{\epsilon_{\alpha}} y_j + D_{ij} \pd{y_j}{\epsilon_{\alpha}} - \pd{Q_i}{\epsilon_{\alpha}} = 0.
\end{equation*}

Also using the solution vector, the sensitivity derivative (iii) is governed by 

\begin{equation*}
    \pd{D_{ij}}{\beta_K} y_j + D_{ij} \pd{y_j}{\beta_K} - \pd{Q_i}{\beta_K} = 0.
\end{equation*}

Using the solution and the sensitivity derivative ii), the derivative (iv) is given by 

\begin{equation*}
    \frac{\partial^2 D_{ij}}{\partial \epsilon_{\gamma} \partial \epsilon_{\alpha}} y_j + \pd{D_{ij}}{\epsilon_{\alpha}} \pd{y_j}{\epsilon_{\gamma}} + \pd{D_{ij}}{\epsilon_{\gamma}} \pd{y_j}{\epsilon_{\alpha}} + D_{ij} \frac{\partial^2 y_j}{\partial \epsilon_{\alpha} \partial \epsilon_{\gamma}} - \frac{\partial^2 Q_i}{\partial \epsilon_{\alpha} \partial \epsilon_{\gamma}} = 0.
\end{equation*}

The mixed sensitivity derivative (v) is computed with 

\begin{equation*}
    \frac{\partial^2 D_{ij}}{\partial \beta_K \partial \epsilon_{\alpha}} y_j + \pd{D_{ij}}{\beta_K} \pd{y_j}{\epsilon_{\alpha}} + \pd{D_{ij}}{\epsilon_{\alpha}} \pd{y_j}{\beta_K} + D_{ij} \frac{\partial^2 y_j}{\partial \epsilon_{\alpha} \partial \beta_K} - \frac{\partial^2 Q_i}{\partial \epsilon_{\alpha} \partial \beta_K}= 0.
\end{equation*}

Finally, the mixed sensitivity derivative which is second-order in the model parameters (vi) is given by 

\begin{equation*}
\begin{aligned}
    \frac{\partial^3 D_{ij}}{\partial \epsilon_{\gamma} \partial \epsilon_{\alpha} \partial \beta_K} y_j + \frac{\partial^2 D_{ij}}{\partial \epsilon_{\gamma} \partial \epsilon_{\alpha}} \pd{y_j}{\beta_K} + \frac{\partial^2 D_{ij}}{\partial \epsilon_{\alpha} \partial \beta_K}\pd{y_j}{\epsilon_{\gamma}} + \pd{D_{ij}}{\epsilon_{\alpha}} \frac{\partial^2 y_j}{\partial \epsilon_{\gamma} \partial \beta_K} + \frac{\partial^2 D_{ij}}{\partial \epsilon_{\gamma} \partial \beta_K} \pd{y_j}{\epsilon_{\alpha}} + \pd{D_{ij}}{\epsilon_{\gamma}} \frac{\partial^2 y_j}{\partial \epsilon_{\alpha} \partial \beta_K} \\ + \pd{D_{ij}}{\beta_K} \frac{\partial^2 y_j}{\partial \epsilon_{\alpha} \partial \epsilon_{\gamma}} + D_{ij}\frac{\partial^3 y_j}{\partial \epsilon_{\alpha} \partial \epsilon_{\gamma} \partial \beta_K} = 0.
\end{aligned}
\end{equation*}



\begin{thebibliography}{10}

\bibitem{brunton_discovering_2016}
Steven~L. Brunton, Joshua~L. Proctor, and J.~Nathan Kutz.
\newblock Discovering governing equations from data by sparse identification of nonlinear dynamical systems.
\newblock {\em Proceedings of the National Academy of Sciences}, 113(15):3932--3937, April 2016.
\newblock Publisher: Proceedings of the National Academy of Sciences.

\bibitem{bunge_inverse_2019}
Mario Bunge.
\newblock Inverse {Problems}.
\newblock {\em Foundations of Science}, 24(3):483--525, September 2019.

\bibitem{chen_hessian_2011}
Pei Chen.
\newblock Hessian {Matrix} vs. {Gauss}–{Newton} {Hessian} {Matrix}.
\newblock {\em SIAM Journal on Numerical Analysis}, 49(4):1417--1435, January 2011.
\newblock Publisher: Society for Industrial and Applied Mathematics.

\bibitem{dashti_bayesian_2015}
Masoumeh Dashti and Andrew~M. Stuart.
\newblock The {Bayesian} {Approach} {To} {Inverse} {Problems}, July 2015.
\newblock arXiv:1302.6989 [math].

\bibitem{hong_global_2020}
Hoon Hong, Alexey Ovchinnikov, Gleb Pogudin, and Chee Yap.
\newblock Global {Identifiability} of {Differential} {Models}.
\newblock {\em Communications on Pure and Applied Mathematics}, 73(9):1831--1879, 2020.
\newblock \_eprint: https://onlinelibrary.wiley.com/doi/pdf/10.1002/cpa.21921.

\bibitem{huan_optimal_2024}
Xun Huan, Jayanth Jagalur, and Youssef Marzouk.
\newblock Optimal experimental design: {Formulations} and computations.
\newblock {\em Acta Numerica}, 33:715--840, July 2024.
\newblock arXiv:2407.16212 [stat].

\bibitem{krampe_optimized_2021}
Valérie Krampe, Pascal Edme, and Hansruedi Maurer.
\newblock Optimized experimental design for seismic full waveform inversion: {A} computationally efficient method including a flexible implementation of acquisition costs.
\newblock {\em Geophysical Prospecting}, 69(1):152--166, 2021.
\newblock \_eprint: https://onlinelibrary.wiley.com/doi/pdf/10.1111/1365-2478.13040.

\bibitem{liu_bayesian_2025}
Tianyi Liu, Xiao Xiao, and Fehmi Cirak.
\newblock Bayesian buckling load optimisation for structures with geometric uncertainties.
\newblock {\em International Journal for Numerical Methods in Engineering}, 126(20):e70111, October 2025.
\newblock arXiv:2501.04553 [math].

\bibitem{maclaren_what_2020}
Oliver~J. Maclaren and Ruanui Nicholson.
\newblock What can be estimated? {Identifiability}, estimability, causal inference and ill-posed inverse problems, July 2020.
\newblock arXiv:1904.02826 [math].

\bibitem{matavel_bayesian-optimized_2025}
Custódio~Efraim Matavel, Andreas Meyer-Aurich, and Hans-Peter Piepho.
\newblock Bayesian-optimized experimental designs for estimating the economic optimum nitrogen rate: a model-averaging approach.
\newblock {\em Agronomy Journal}, 117(3):e70087, 2025.
\newblock \_eprint: https://acsess.onlinelibrary.wiley.com/doi/pdf/10.1002/agj2.70087.

\bibitem{mentre_sparse-sampling_1995}
France Mentré, Pascale Burtin, Yann Merlé, Joost van Bree, Alain Mallet, and Jean-Louis Steimer.
\newblock Sparse-{Sampling} {Optimal} {Designs} in {Pharmacokinetics} and {Toxicokinetics}*.
\newblock {\em Drug Information Journal}, 29(3):997--1019, July 1995.
\newblock Publisher: SAGE Publications.

\bibitem{mercier_designing_2025}
Arnaud Mercier, Christian Boehm, and Hansruedi Maurer.
\newblock Designing full waveform inverse problems: a combined data and model approach.
\newblock {\em Geophysical Journal International}, 241(3):1479--1494, June 2025.

\bibitem{miao_identifiability_2011}
Hongyu Miao, Xiaohua Xia, Alan~S. Perelson, and Hulin Wu.
\newblock On {Identifiability} of {Nonlinear} {ODE} {Models} and {Applications} in {Viral} {Dynamics}.
\newblock {\em SIAM Review}, 53(1):3--39, January 2011.
\newblock Publisher: Society for Industrial and Applied Mathematics.

\bibitem{pukelsheim_optimal_1993}
Friedrich Pukelsheim.
\newblock {\em Optimal {Design} of {Experiments}}.
\newblock SIAM Publications, 1993.

\bibitem{raissi_physics-informed_2019}
M.~Raissi, P.~Perdikaris, and G.E. Karniadakis.
\newblock Physics-informed neural networks: {A} deep learning framework for solving forward and inverse problems involving nonlinear partial differential equations.
\newblock {\em Journal of Computational Physics}, 378:686--707, February 2019.

\bibitem{ricciardi_bayesian_2023}
Denielle Ricciardi, Tom Seidl, Brian Lester, Amanda Jones, and Elizabeth Jones.
\newblock Bayesian {Optimal} {Experimental} {Design} for {Constitutive} {Model} {Calibration}, October 2023.
\newblock arXiv:2308.10702 [cs].

\bibitem{rowan_extending_2025}
Conor Rowan.
\newblock Extending the explicit constraint force method to inverse problems, December 2025.
\newblock arXiv:2512.14877 [math].

\bibitem{rowan_okay_2025}
Conor Rowan.
\newblock {\em Okay, this one is made up.}
\newblock 2025.

\bibitem{rowan_physics-informed_2025}
Conor Rowan, Kurt Maute, and Alireza Doostan.
\newblock Physics-informed solution reconstruction in elasticity and heat transfer using the explicit constraint force method, May 2025.
\newblock arXiv:2505.04875 [cs].

\bibitem{senig_investigation_2025}
James Davis~A. Senig and John~F. Maddox.
\newblock An {Investigation} {Into} the {Effective} {Gaseous} {Thermal} {Conductivity} of {Fibrous} {Insulation} {Materials}.
\newblock In {\em {AIAA} {AVIATION} {FORUM} {AND} {ASCEND} 2024}. American Institute of Aeronautics and Astronautics, 2025.
\newblock \_eprint: https://arc.aiaa.org/doi/pdf/10.2514/6.2024-4030.

\bibitem{sun_local_2022}
Jiguang Sun.
\newblock Local estimators and {Bayesian} inverse problems with non-unique solutions.
\newblock {\em Applied Mathematics Letters}, 132:108149, October 2022.

\bibitem{tao_digital_2019}
Fei Tao, He~Zhang, Ang Liu, and A.~Y.~C. Nee.
\newblock Digital {Twin} in {Industry}: {State}-of-the-{Art}.
\newblock {\em IEEE Transactions on Industrial Informatics}, 15(4):2405--2415, April 2019.
\newblock Conference Name: IEEE Transactions on Industrial Informatics.

\bibitem{vogel_computational_2002}
Curtis~R. Vogel.
\newblock {\em Computational {Methods} for {Inverse} {Problems}}.
\newblock Frontiers in {Applied} {Mathematics}. Society for Industrial and Applied Mathematics, January 2002.

\end{thebibliography}

\end{document}